\documentclass[a4paper]{article}

\usepackage[english]{babel}
\usepackage[utf8x]{inputenc}
\usepackage[T1]{fontenc}

\usepackage[margin=2cm, head=1cm]{geometry}

\usepackage{amsmath,amssymb,graphicx,url,natbib}
\usepackage{subfigure}
\usepackage{lscape}
\usepackage{adjustbox}
\usepackage[colorinlistoftodos]{todonotes}
\usepackage[colorlinks=true, allcolors=blue]{hyperref}

\newtheorem{theorem}{Theorem}

\begin{document}

\title{Controlling IL-7 injections in HIV-infected patients
}


\author{Chlo\'{e} Pasin$^{1,2,3,4}$, Fran\c{c}ois Dufour$^{1,2,5}$, Laura Villain$^{1,2,3,4}$, Huilong Zhang$^{1,2}$, Rodolphe Thi\'{e}baut$^{1,2,3,4,\dagger}$}

\date{January 2018}

\maketitle

\noindent
$^{1}$ Univ Bordeaux, France \\
$^{2}$ INRIA Bordeaux Sud Ouest, Talence, France \\
$^{3}$ ISPED, Centre INSERM U1219, Bordeaux, France \\
$^{4}$ Vaccine Research Institute, VRI, H\^opital Henri Mordor, Cr\'eteil France \\
$^{5}$ Bordeaux INP, IMB, Bordeaux, France \\
$^{\dagger}$ email: rodolphe.thiebaut@u-bordeaux.fr

\begin{abstract}
Immune interventions consisting in repeated injections are broadly used as they are thought to improve the quantity and the quality of the immune response. However, they also raised several questions that remains unanswered, in particular the number of injections to make or the delay to respect between different injections to achieve this goal. Practical and financial considerations add constraints to these questions, especially in the framework of human studies. We specifically focus here on the use of interleukine-7 (IL-7) injections in HIV-infected patients under antiretroviral treatment, but still unable to restore normal levels of CD4$^+$ T lymphocytes. Clinical trials have already shown that repeated cycles of injections of IL-7 could help maintaining CD4$^+$ T lymphocytes levels over the limit of 500 cells/$\mu$L, by affecting proliferation and survival of CD4$^+$ T cells. We then aim at answering the question : how to maintain a patient’s level of CD4$^+$ T lymphocytes by using a minimum number of injections (ie optimizing the strategy of injections) ? Based on mechanistic models that were previously developed for the dynamics of CD4+ T lymphocytes in this context, we model the process by a piecewise deterministic Markov model. We then address the question by using some recently established theory on impulse control problem in order to develop a numerical tool determining the optimal strategy. Results are obtained on a reduced model, as a proof of concept : the method allows to define an optimal strategy for a given patient. This method could be applied to optimize injections schedules in clinical trials.

\end{abstract}

\section{Introduction}
\label{intro}
The infection by the Human Immunodeficiency Virus (HIV) compromises the immune system functions, mainly because of the depletion of CD4$^+$ T lymphocytes. Combined antiretroviral (cART) therapy has led to a spectacular improvement of patients' survival by controlling virus replication and consequently restoring the immune system functions. However, some patients fail at reconstituting their immune system and recovering normal CD4$^+$ T cell levels, especially when they start antiretroviral treatment late (\cite{lange2003immune}), that is when already reaching low CD4+ T cell count. Immune therapy has been considered as a complement to cART to help immune restoration. In particular, interleukin-7 (IL-7), a cytokine produced by non-marrow-derived stromal and epithelial cells, is thought to improve thymic production (\cite{mackall20017}; \cite{okamoto2002effects}) and cell survival (\cite{tan20017}; \cite{vella1998cytokine}, \cite{leone2010increased}). The safety and beneficial effect of injections of exogenous IL-7 was first shown in phase I trials (\cite{sereti20097}; \cite{levy2009enhanced}) and observational studies (\cite{camargo2009responsiveness}). Then, phase I/II human clinical trials (INSPIRE 1, 2 and 3 studies) have evaluated the effect of repeated cycles of three IL-7 injections and showed that this therapy helped maintaining HIV infected patients with CD4$^+$ T cells levels above 500 cells/$\mu$L (\cite{levy2012effects}), a level associated with a nearly healthy clinical status (\cite{lewden2007hiv}).\\

The dynamics of CD4$^+$ T lymphocytes following IL-7 injections can be fitted by mechanistic models based on ordinary differential equations (ODEs), with compartments corresponding to different populations of CD4$^+$ T lymphocytes and biological parameters characterizing these populations. Hence it was possible to quantify the effect of repeated cycles of IL-7 on CD4$^+$ T lymphocytes on specific parameters. Previous work using data from clinical trials (INSPIRE studies) has shown that IL-7 enhances both proliferation and survival of CD4$^+$ T lymphocytes (\cite{thiebaut2014quantifying}). Moreover, a differential effect of the injections within a given cycle have been found, the third injection of a cycle appearing to have a weaker effect on proliferation than the first ones (\cite{jarne2017modeling}). \\

In addition to providing insight into the most important mechanism of the effect of exogenous IL-7, the models have shown a very good predictive capacity (\cite{thiebaut2014quantifying}; \cite{jarne2017modeling}). Hence, the next step was the determination of the best protocol of injections. A first approach, realized in \cite{jarne2017modeling}, consisted in simulating and comparing three protocols to the regular one. In all four protocols, CD4 counts were measured every three months, and a new cycle was administered when the CD4 numbers were below 550 cells/$\mu$L. Simulated protocols consisted in either a first three-injections cycle followed by repeated two-injections cycles, or a first three-injections cycle followed by repeated one-injection cycles, or two-injections cycles. Comparison was based on three criteria : number of injections received, mean CD4 count and time spent below 500 cells/$\mu$L over a four-year period. Results showed that cycles of two injections could be sufficient to maintain CD4 levels, while using less injections than in the clinical protocol. However, even if these results allow to consider reducing the number of injections in clinical protocols, there are still the constraints due the delay between visits, as follow-up is based on measurements every three months (independently of the patient). Indeed, some "not too low responders" could afford coming back later than three months after the last visit whereas some "low responders" would need more repeated cycles or more injections by cycle: in the simulated protocols, the so-called "very low responders" still had to deal with  CD4 levels under 500 cells/$\mu$L during a meantime of 400 days over the four-year period of study. Thus, individualized protocols could help achieving the maintenance of the patient's CD4$^+$ T lymphocytes levels over a given threshold by using different patient-dependent timing of injections and doses. The possibility of conducting the lightest intervention for every patients could be very important for the development of IL-7 in HIV infected patients especially for further large clinical trials. Optimization of schedule and doses is an up-to-date question when working on protocol of injections. In their review on mathematical modeling for immunology, \cite{eftimie2016mathematical} emphasize the need for complex optimal control approaches coupled with immunology experiments, in order to improve clinical interventions. Some work using optimal control has already been done in cancer immunotherapy as in \cite{castiglione2006optimal}, \cite{cappuccio2007determination},  \cite{castiglione2007cancer} and \cite{pappalardo2010vaccine} to determine the optimal schedule of administration of dendritic cells vaccine, or in drug dosage for in vitro fertilization, as in \cite{yenkie2014optimal}. In the HIV field, \cite{stengel2008mutation},  \cite{yang2013pulse} and \cite{croicu2015short} were interested in determining an optimal HIV treatment to minimize the viral load and maximize the number of uninfected CD4$^+$ T cells taking into account the occurrence of viral mutations. Basically, there are two kinds of techniques that can be used to solve optimal control problems: methods involving Pontryagin's maximum principle and dynamic programming approaches. Pontryagin's maximum principle provides a set of necessary conditions (sometimes sufficient) that need to be satisfied by an optimal control and the optimal problem can be solved using the associated Hamiltonian system. Pontryagin's maximum principle has been applied to a number of biological problems of the form $\frac{dx(t)}{dt}=f(x(t),u(t))$, where the solution to the ordinary differential equation depends on the dynamics of the control function $u(t)$. For example, this is the case in \cite{castiglione2006optimal}, \cite{cappuccio2007determination},  \cite{castiglione2007cancer} and \cite{pappalardo2010vaccine}, where the dynamics of dendritics cells are modeled as a death term plus an input term induced by the vaccine injection. To account for the uncertainty and measurement that can induce errors, \cite{stengel2004stochastic} worked on a stochastic therapeutic protocol and adapted the method by minimizing the expected value of the cost subject to a stochastic constraint.  However, in our case, the model is a piecewise deterministic Markov model (PDMP), where dynamics of IL-7 are unknown and thus not modeled. Addressing the objective of spending the least time possible under the threshold of 500 cells/$\mu$L by using repeated injections of IL-7 corresponds in a more formal way to determining actions (injection or not and choice of dose) at given time points over an horizon of time: this can be treated as a problem of impulse control in the optimal control theory. To the best of our knowledge, there is no maximum principle solving this kind of problem. We will then focus on a dynamic programming method, as developed in \cite{costa2016constrained}. In a formal mathematical framework, we addressed the question of optimizing the schedule of IL-7 injections for a given patient by a two-steps method : determining an adapted mathematical model for the process, and developing a numerical method to determine an optimal strategy of IL-7 injections for a given patient.\\

As described in \cite{davis1984piecewise}, most of the continuous-time stochastic problems of applied probability (thus including those modeling biological processes) consist of some combination of diffusion, deterministic motion and/or random jumps. Ordinary differential equations can be included in the class of deterministic motion with random jumps. In our particular framework of modeling cell dynamics after IL-7 injections, jumps correspond to the change of some parameters value. This can be easily and naturally modeled by the largely studied class of Piecewise Deterministic Markov Processes (PDMPs). A non-controlled version of this model can be described by iteration as follows : from a point in the state space, the process follows a deterministic trajectory determined by the flow, until a jump occurs. This jump happens either spontaneously in a random manner, or when the flow hits the boundary of the state space. After the jump, the system restarts from a new point determined by the transition measure of the process. We will show in this article how to model the dynamics of the CD4$^+$ T cells in HIV-infected patients following IL-7 injections using a PDMP. \\

According to the problem studied in \cite{costa2016constrained}, impulse control consists in possible actions only when the process reaches its boundary. This will constitute our framework : the decision maker has the possibility to inject IL-7 when the number of CD4$^+$ T lymphocytes reaches a given level or when a certain amount of time has passed since the last injection, which will correspond to the boundary of our state space. From a mathematical point of view, we note $X$ the state space, an open subset of $\mathbb{R}^d$, $d \in \mathbb{N}$, and $\partial X$ its boundary. The flow associated to the process is $\phi(x,t) : \mathbb{R}^d \times \mathbb{R} \mapsto \mathbb{R}^d$.  The active boundary is defined as $\Xi = \{x \in \partial X : x=\phi(y,t)\}$ for some $y \in X$ and $t \in \mathbb{R}_+^*$. We will then denote $\overline{X}=X \cup \Xi$ and for $x \in \overline{X}$, we can define $t^*(x)=\inf \{t \in \mathbb{R}_+ : \phi(x,t) \in \Xi\}$. The controlled jump intensity $\eta$ is a $\mathbb{R}_+$-valued measurable function and determines the law of the stochastic jumps. When the process, i.e. the trajectory of CD4$^+$ T lymphocytes, reaches $\Xi$, the decision-maker can act on the patient by injecting IL-7. The action is thus the possible dose injected. This will lead to a jump in some parameters value, and the process will restart from a new point defined by the transition measure $Q(.|\phi(x_0,\tau),d)$, depending on the dose and the position of the state before the jump $\phi(x_0,\tau)$. Therefore, the impulse control problem consists in determining the optimal scheme of injections and their associated dose according to a given optimality criterion, based on the cost function $\mathit{C}$: in our case, this cost function will depend on the number of injections realized and the time spent with the CD4$^+$ T lymphocytes levels under the threshold of 500 cells/$\mu$L. Both quantities need to be minimized, in order to maintain the patient in good health by injecting the least possible. The set of all realized injections over a given horizon constitute the policy. In a more formal way, a strategy $u$ of the decision maker is a sequence $u = \{u_n\}_{n \in \mathbb{N}}$ of functions $u_n : X \mapsto A$ giving the action to realize at punctual time points $t_n \geq 0$ when the system is in state $x \in X$. The impulse control problem aims at finding a policy $u$ minimizing the discounted cost optimality criterion defined as: 
\begin{equation*}
\mathit{V}(u,x_0)=\mathbb{E}_{x_0}^{u}\Big[ \sum_n e^{-\alpha t_n} \mathit{C}(x_n,a_n)\Big]
\end{equation*}
where $\mathbb{E}_{x_0}^{u}$ is the expectation operator under policy $u$ and initial state $x_0$ and $0<\alpha<1$ a discount factor. In our case, the impulse control problem consists then in determining the patient-specific schedule of injections and their dose to optimize the patient's CD4$^+$ T lymphocyte numbers by using a minimum number of injections. \\

As emphasized by the authors of \cite{dufour2015numerical}, the development of computation methods for the control of PDMPs has been limited, and at the moment there is no general method allowing the numerical resolution of optimal control on PDMPs (and in particular impulse control). This corresponds to a real challenge. We propose in this work a numerical method based on the results developed in \cite{costa2016constrained}. In this paper, the authors studied the existence of a solution of the Bellman-Hamilton-Jacobi equation by showing that the value function is the limit of a sequence of functions given by iteration of an integro-differential operator. This construction leads to a natural method for the computation of the optimal cost and the determination of an optimal strategy of injections. In particular, we applied our numerical tool to the case of a simplified biological model as a proof of concept and succeeded in providing an optimal strategy for two given patients with different profiles. This paper is organized as follows : section 2 presents the mathematical modeling of the process, including data and design of INSPIRE studies, as well as mechanistic model and finally the associated PDMP. Section 3 focuses on the optimal control problem, by reminding the main theoretical results from \cite{costa2016constrained} and adapting them to the IL-7 study. Section 4 presents some numerical aspects of the dynamic programming work, necessary to determine the optimal cost function and strategy for a given patient. Results are presented in section 5 and discussion is done in section 6.

\section{Mathematical modeling}
\label{sec:1}
\subsection{Material}
\label{sec:1.1}
Our work is based on three phase I/II multicenter studies assessing the effect of a purified glycosylated recombinant human Interleukin 7 (IL-7) treatment for immune restoration in HIV infected patients under treatment: INSPIRE (\cite{levy2012effects}), INSPIRE 2 and INSPIRE 3 (\cite{thiebaut2016repeated}). A total of 128 HIV-infected patients under antiretroviral therapy with CD4$^+$ T cell count between 100-400 cells/$\mu$L and undetectable viral load for at least 6 months were included among the three studies from the time of the first injection. IL-7 was administered in cycles of weekly injections, with a "complete cycle" defined as three weekly injections. In INSPIRE, all 21 patients received complete cycles of IL-7 at different weight-dependent doses: 10, 20 and 30 $\mu$g/kg. In INSPIRE 2 and INSPIRE 3, 23 and 84 patients (respectively) received repeated (and sometimes incomplete) cycles of IL-7 at dose 20 $\mu$g/kg. Repeated visits and follow-up once every 3 months after the first cycle allowed to measure biomarkers levels in patients, in particular total CD4$^+$ T cell counts and number of proliferating CD4$^+$ T cells through Ki67 marker. At every visit, a new cycle of injections was administered if the patient's CD4$^+$ T cell level was under 550/$\mu$L, in order to globally maintain the levels above 500 cells/$\mu$L. The total duration of the studies was 12, 24 and 21 months for INSPIRE, INSPIRE 2 and INSPIRE 3, respectively. 

\subsection{Mechanistic model}
\label{sec:1.2}
\begin{figure}
\centering
\includegraphics[width=0.5\textwidth]{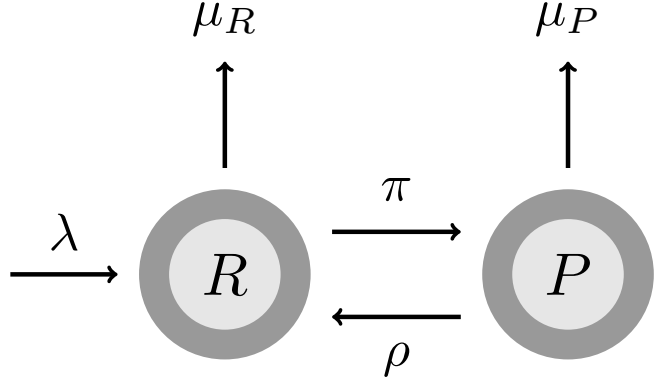}
\caption{Mechanistic model for the dynamics of CD4$^+$ T lymphocytes}
\label{mechamodel}      
\end{figure}

The dynamics of CD4$^+$ T lymphocytes were largely described in \cite{thiebaut2014quantifying} and \cite{jarne2017modeling} by using several mechanistic models. We will focus here on the following model, described in figure \ref{mechamodel}, 
including two populations of cells, non-proliferating (or resting, R) and proliferating (P). Resting cells are produced by thymic output at rate $\lambda$, become proliferating cells at rate $\pi$ and die at rate $\mu_R$. Proliferating cells die at rate $\mu_P$ and can also divide and produce two non-proliferating cells at rate $\rho$. The system of differential equations is written as follows :
\begin{equation}
\label{odesystem}
  \left\{
    \begin{aligned}
      \frac{dR}{dt}&=\lambda - \mu_R R - \pi R + 2 \rho P \\
      \frac{dP}{dt}&= -\mu_P P - \rho P + \pi R\\
    \end{aligned}
  \right.
\end{equation}
We assume the system is at equilibrium at $t=0$, before the study begins and any injection is administered. IL-7 injections are realized through cycles containing up to three injections with seven days elapsed between each injection. Parameters estimation was performed using a population approach. Mixed-effect models including intercept, random and fixed effects, were used on log-transformed parameters, in order to both obtain an estimation across population and account for between-individuals variability. In the controlled framework, the decision-maker can decide to inject IL-7 to a patient at a given dose $d$, and this will affect the value of the proliferation rate $\pi$. Each injection denoted by $n \in \{1,2,3\}$ of a cycle has a different effect on the value of $\pi$ for patient $i$, defined as follows :
\begin{equation}
\label{injeffpi}
\tilde{\pi}^i = \tilde{\pi}_0 + \beta_{\pi}^{(n)} d_i^{0.25} \mathbf{1}_{\{t \in [t^i_{inj},t^i_{inj}+\tau^i]\}} 
\end{equation}
with $\tilde{\pi}=log(\pi)$; $\beta_{\pi}=[\beta_{\pi}^{(1)},\beta_{\pi}^{(2)},\beta_{\pi}^{(3)}]$ is the vector of effect of each injection of a single cycle; $d$ is the injected dose; $t_{inj}$ is the time (in days) at which IL-7 is injected and $\tau$ the length of effect of the injection (in number of days), considered equal to 7 in previous models \cite{jarne2017modeling}. Estimation of parameters showed that effect of successive injections on the proliferation rate decreases within a cycle, and in particular the third injection seems to have a much weaker effect (thus $\beta_{\pi}^{(3)} < \beta_{\pi}^{(2)} < \beta_{\pi}^{(1)}$) (\cite{jarne2017modeling}). 

\subsection{Mathematical model : piecewise deterministic Markov process}
\label{sec:1.3}
As described in the introduction, ODEs-based mechanistic models can be included into the broader class of PDMPs. A PDMP is characterized by a state space in which it evolves, a flow, a jump intensity and a measure of transition. In this section, we present  the PDMP associated to the biological process described in section \ref{sec:1.2}. The PDMP modeling the dynamics of CD4$^+$ T lymphocytes of a given patient is defined using six variables: the state vector is then denoted by $x=(\gamma,n,\sigma,\theta,p,r)$, with $\gamma$ determining the value of parameter $\pi$ when combined with $n$, the number of injections realized in the ongoing cycle. If $d=[d_0,d_1,..,d_{m_d}]$ is the vector of all possibles doses (with $d_0=0$), then $\gamma \in \{1..m_d+1\}$. Injecting dose $d_k$ at the $n$-th injection of a cycle gives the following : $\gamma(d_k)=k+1$ and $\pi=\pi_0+\beta_{\pi_n}d(\gamma)^{0.25}$. The two variables $\sigma$ and $\theta$ are time variables, discretized with step of one day. In particular, $\sigma$ corresponds to the number of days since the last injection and $\theta$ to the running time ($\theta=1$ at the first injection of the first cycle). Finally, variables $p$ and $r$ are values of compartments P and R solutions of system \ref{odesystem} with parameter $\pi$ defined by $\gamma$ and $n$. We suppose the patient is followed until an horizon of time $T_h$, then the state space is $X=\tilde{X} \cup \Delta$ with 
\begin{equation*}
    \begin{array}{rl}
        \tilde{X} =& \{1..m_d+1\} \times \{1..n_{inj}\} \times \{0..T_h-7(n_{inj}-1)\} \times \{0..T_h\}\\
         &\times \{p_{min}..p_{max}\} \times \{r_{min}..r_{max}\}\\         
    \end{array}
\end{equation*}
and $\Delta$ is an absorbing state representing the end of the study, at $t=T_h$ : $\Delta=(0,0,0,T_h,0,0)$. For $x=(\gamma,n,\sigma,\theta,p,r) \in X$, the flow is defined as : 
\begin{itemize}
\item $\phi(x,t)=(\gamma,n,\sigma+t,\theta+t,p,r)$ if $\theta \leq 1$
\item $\phi(x,t)=(\gamma,n,\sigma+t,\theta+t,P(t,\gamma,n),R(t,\gamma,n))$ if $\theta \in [1,T_h-1]$, with $P(t,\gamma,n)$ and $R(t,\gamma,n)$ solutions of system \ref{odesystem} with initial conditions $p$ and $r$ and $\pi$ determined with $\gamma$ and $n$
\item $\phi(\Delta,t)=\Delta $
\end{itemize}
Moreover, even if the deterministic mechanistic model allowed good fits for the data, we make the hypothesis that the process undergoes some stochasticity: in particular, as the value of parameter $\pi$ is modified by an injection of IL-7 during some days, we suppose that this modification can last randomly up to 7 days after the injection. Stochastic jumps can then occur with intensity $\eta$ such that for $x \in \tilde{X}$, $\eta(x)=\eta \mathbf{1}_{\{\gamma >1\}} $ with $\eta$ a given value and $\eta(\Delta)=0$. It means that if we consider the modification of $\pi$ value after an IL-7 injection through equation \ref{injeffpi}, $\tau$ follows there a random exponential law of parameter $\eta$. We define the constant $K=\eta$ such that $\eta(x) \leq K$ for every $x\in \tilde{X}$.\\
As the main objective of the injections is to maintain the CD4+ T cell level over 500 cells/$\mu$L, we define the deterministic boundary of the state space when this value is reached: this will allow a possible new IL-7 injection that should increase the CD4+ T cells count. In order to account for the clinical constraint, we assumed a minimum time $\sigma_{min}$ is observed between the beginning of two consecutive cycles, even if the number of CD4 falls below the threshold of 500. During cycles, the deterministic boundary reflects the seven days delay between injections. In a more formal way, the boundary can actually be reached in five different situations described in the following : 
\begin{itemize}
\item for a technical reason due to the mathematical modeling which cannot account for an impulse action at $t=0$, we define a first artificial boundary when the study begins, at $\theta=1$ : $\Xi_1=\{x : \theta=1\}$. This allows a cycle of injections to begin at $\theta=1$. We suppose the studied patient is already included in the clinical study: it means that its biological parameters are known, and her/his CD4$^+$ T cell count at $t=0$ as well (either because she/he is at equilibrium, and the values are known from biological parameters, or because some measures have been realized at this time).
\item we also need to determine a time corresponding to the end of the study, so we define a boundary when the time reaches the horizon $T_h$ :$\Xi_2=\{x : \theta \geq T_h\}$
\item another boundary is reached when the patient is undergoing a cycle of injections and seven days have passed since the last injection : $\Xi_3=\{x : n<n_{inj}, \sigma=7, \theta < T_h \}$
\item we also consider when at least one cycle was already achieved and the count of cells is equal to or below the threshold of 500 cells$/\mu L$. We also assumed a minimum time $\sigma_{min}$ is observed between the beginning of two consecutive cycles : $\Xi_4 = \{ x : p+r \leq 500, n=n_{inj}, \sigma \geq \sigma_{min}, \theta < T_h\}$
\item finally, an artificial boundary is created when $\pi$ has not returned to its baseline value seven days after the last injection of a cycle : $\Xi_5=\{x : \gamma >1, n=n_{inj}, \sigma=7, \theta<T_h \}$
\end{itemize}
We define then the active boundary as $\Xi=\Xi_1 \cup \Xi_2 \cup \Xi_3 \cup \Xi_4 \cup \Xi_5$.
In this process, actions (IL-7 injections) can only be realized when the process hits the active boundary. We model the possibility of not doing an injection in a given cycle by using a fictive dose $d_0$ equal to zero. When beginning a new cycle of injections, the first injection needs to be positive though. The possible action that can be chosen by the decision-maker depends on the boundary reached. Therefore for every $x \in \Xi$ : 
\begin{equation*}
A(x) = \left\{
    \begin{array}{ll}
        \{d_1,..d_{m_d}\} & \hbox{\ if\ } x \in \Xi_1 \cup \Xi_4\\
        \{0,d_1,..d_{m_d}\} &\hbox{\ if\ } x \in \Xi_3 \\
         \emptyset & \hbox{\ if\ } x \in \Xi_2 \cup \Xi_5 
    \end{array}
\right.
\end{equation*}
Now we can define the transition measure (or Kernel), which determines the new point from which the process restarts after a jump. It depends on the injected dose only when the boundary of the process is reached. All possible situations are the following : 
\begin{itemize}
\item when the flow hits $\Xi_1$, the study begins with administration of a cycle of injections. $\gamma$ takes the value corresponding to the chosen dose. $(p,r)=(P_c,R_c)$, known values from either equilibrium or biological measures made on the patient before the beginning of the study
\item when the flow hits $\Xi_2$, the study is over and nothing happens from absorbing state $\Delta$
\item when the flow hits $\Xi_3$, a new injection is administered to the patient. $\gamma$ takes the value corresponding to the chosen dose $\gamma(d)$, $n$ increases by one, $\sigma$ goes back to 0
\item when the flow hits $\Xi_4$, a new cycle of injections begins. $\gamma$ takes the value corresponding to the chosen dose, $n$ goes back to 1, $\sigma$ goes back to 0
\item when the flow hits $\Xi_5$, there is no injection. $\gamma$ goes back to 1
\item in case of spontaneous jump, there is no injection and $\gamma$ goes back to 1
\end{itemize}
In a formal way, the Kernel Q is written : 
\begin{equation*}
\begin{array}{rll}
Q(dy|x,d) &= &
        \delta_{(\gamma(d),1,0,1,P_c,R_c)}(dy)\mathbf{1}_{\{x \in \Xi_1\}} + \delta_{\Delta}(dy)\mathbf{1}_{\{x \in \Xi_2\}}\\
       &+& \delta_{(\gamma(d),n+1,0,\theta,p,r)}(dy)\mathbf{1}_{\{x \in \Xi_3\}} +\delta_{(\gamma(d),1,0,\theta,p,r)}(dy)\mathbf{1}_{\{x \in \Xi_4\}}   \\
   & +& \delta_{(1,n,\sigma,\theta,p,r)}(dy)\mathbf{1}_{\{x \in \Xi_5\}} + \delta_{(1,n,\sigma,\theta,p,r)}(dy)\mathbf{1}_{\{x \in \tilde{X}\}}
    \end{array}
\end{equation*}

As explained in the introduction, the objective of this work is to maintain the CD4$^+$ T cell count above 500 cells/$\mu$L while using a minimum number of injections. The determination of the optimal strategy relies on the optimization of a criterion function, divided in two costs. First, the gradual cost penalizes the trajectory of the process through the time spent under the threshold after the beginning of the first cycle. This time is considered in months, approximately, as it is computed as the number of days divided by 30. Then for $x=(\gamma,n,\sigma,\theta,p,r) \in \tilde{X}$ : 
\begin{equation*}
C^g(x) =  \frac{1}{30}\mathbf{1}_{\{p+r<500\}}\mathbf{1}_{\{\theta \geq 1\}} 
\end{equation*}
Then, the cost associated to an impulsive action penalizes the fact of injecting IL-7 to the patient : 
\begin{equation*}
C^i(x,d)= \mathbf{1}_{\{x \in \Xi_1 \cup \Xi_4\}} + \mathbf{1}_{\{d \neq 0\}}\mathbf{1}_{\{x \in \Xi_3\}}
\end{equation*}
After the horizon, the cost is null, as $C^i(\Delta)=C^g(\Delta)=0$. 

\section{Optimal control}
\label{sec:2}
In this section, we will first remind the main theoretical results obtained in \cite{costa2016constrained}, then we will transpose these results to our particular context.

\subsection{Main theoretical results}
\label{sec:2.1}
The objective of this section is to adapt some results obtained in \cite{costa2016constrained} to our specific context. We follow closely their notation. As defined in the introduction, an admissible control strategy (with only impulse actions) is a sequence $u=(u_n)_{n \in \mathbb{N}}$ of IL-7 doses to be injected at some points of the state space. The set of all admissible strategies is noted $\mathcal{U}$. According to section 2.2 in \cite{costa2016constrained} there exists a continuous-time stochastic process $\xi$ defined on probabilistic space using characteristics $\phi$,$\eta$ and $Q$ depending on the action given by $u$, such that $\xi_t$, $t \in \mathbb{R}_+$ corresponds to the state of the variables at time $t$. The cost of an admissible strategy $u \in \mathcal{U}$ depends on the gradual cost on the trajectory of the process $\xi$, $C^g$, and the cost related to an injection, $C^i$, as defined in section \ref{sec:1.3}:
\begin{equation*}
\begin{array}{rll}
\mathit{V}(u,x_{0}) &= &\mathbb{E}^{u}_{x_{0}} \Bigg[  \int_{]0,+\infty[} e^{-\alpha s}  C^{g}(s) ds \Bigg] \\
&&+ \mathbb{E}^{u}_{x_{0}} \Bigg[  \int_{]0,+\infty[} e^{-\alpha s} I_{\{\xi_{s-}\in\mathbf{\Xi}\}}\int_{\mathbf{A}(\xi_{s-})} C^{i}(\xi_{s-},a) u(da|s) \mu(ds) \Bigg]
\label{Def-cost-j}
\end{array}
\end{equation*}
with $\alpha > 0$ the discount factor and where $\mu$ is a measure that counts the number of jumps in the process. The theorem allowing to determine the optimal cost and providing an optimal strategy is the following:
\begin{theorem} 
\label{maintheorem}
Suppose assumptions A, B and C from section 3.2 in \cite{costa2016constrained} are verified. We define the sequence of functions $\{W_q\}_{q \in \mathbb{N}}$ for any $x\in \overline{\mathbf{X}}$ as follows: 
\begin{equation}
\left\{
    \begin{array}{rll}
        W_{q+1}(x)& =&\mathfrak{B}W_q(x) \hbox{\ for\ q} \in \mathbb{N}\\
        W_0(x) &=&-K_{A}\mathbf{1}_{A_{\varepsilon_1}}(x)-(K_{A}+K_{B})\mathbf{1}_{A_{\varepsilon_1}^{c}}(x) 
    \end{array}
\right.
\label{seqW}
\end{equation}
with constants $K_A$ and $K_B$ defined as in section 5 of \cite{costa2016constrained}, $A_{\varepsilon_1}=\{ x\in\mathbf{X}:t^*(x)>\varepsilon_1\}$ and 
\begin{equation}\label{Def-BV}
\mathfrak{B}V(y)=
\int_{[0,t^{*}(y)[} e^{-(K+\alpha)t} \mathfrak{R}V(\phi(y,t)) dt + e^{-(K+\alpha)t^{*}(y)} \mathfrak{T}V(\phi(y,t^{*}(y)))
\end{equation}
with real-value functions $\mathfrak{R}V$ and $\mathfrak{T}V$ defined for any $V$ respectively on X and $\Xi$ : 
\begin{equation*}
\mathfrak{R}V(x)  =  C^{g}(x) + qV(x)+\eta V(x)
\end{equation*}
\begin{equation*}
\mathfrak{T}V(z) =  \inf_{d\in \mathbf{A}(z)} \Big\{ C^{i}(z,d) + QV(z,b)\Big\}
\end{equation*}
$q$ being the signed kernel, which computes the difference between the states before and after the spontaneous jump. For $x \in X$, it is defined with : 
\begin{equation*}
q(dy|x) = \eta(x)[Q(dy|x)-\delta_x(dy)] 
\end{equation*}
The sequence of functions $\{W_q\}_{q \in \mathbb{N}}$ converges to a function $W$ defined on the state space and such that : 
\newline
i) $W(x_0)=\inf_{u \in \mathcal{U}} \mathcal{V}(u,x_0)$
\newline
ii) there is a measurable mapping $\widehat{\varphi}: \Xi \to \mathbf{A}$ such that $\widehat{\varphi}(z)\in \mathbf{A}(z)$ for any $z \in \Xi$ and satisfying 
\begin{equation}
\label{action_eq}
C^{i}(z,\widehat{\varphi}(z)) + QW(z,\widehat{\varphi}(z))=\inf_{d\in \mathbf{A}(z)} \Big\{ C^{i}(z,d) + QW(z,b)\Big\}.
\end{equation}
\end{theorem}
This theorem allows to determine the optimal cost and an optimal injection strategy, consisting in choosing the optimal action $\widehat{\varphi}(z)$ for every point $z \in \Xi$ reached on the trajectory of the process. Indeed, the iteration of the sequence $\{W_q\}_{q \in \mathbb{N}}$ defined by equation \ref{seqW} can be realized by numerically approximating the operator $\mathfrak{B}$ defined in equation \ref{Def-BV}. This will give an approximation of the function $W$, and in particular of $W(x_0)$, corresponding to the optimal cost. Moreover, to obtain an optimal strategy, the process is the following : we simulate a trajectory from $x_0$, then when a boundary is reached, the chosen action corresponds to the one minimizing the criterion $C^{i}(z,d) + QW(z,b)$, as given by equation \ref{action_eq}. 

\subsection{Application}
\label{sec:2.2}
The process describing the effect of IL-7 on CD4$^+$ T lymphocytes dynamics is now well defined by its characteristics $\phi$, $\eta$ and $Q$, boundaries and possible actions in section \ref{sec:1.3}. Moreover, both gradual cost on the trajectory and impulse cost were defined in that section. We will quickly describe in this part how to apply the results from theorem \ref{maintheorem} for our specific problem, i.e. determining the function $\mathfrak{B}$ needed for the computation of the optimal strategy. For a more detailed and formal computation, we refer the reader to appendix \ref{computeB}. We need to compute:
\begin{equation*}
\mathfrak{B}V(y)=
\int_{[0,t^{*}(y)[} e^{-(K+\alpha)t} \mathfrak{R}V(\phi(y,t)) dt + e^{-(K+\alpha)t^{*}(y)} \mathfrak{T}V(\phi(y,t^{*}(y)))
\end{equation*}
We define 
\begin{equation}
\label{defG}
G(V,y)=\int_{[0,t^*(y)[} e^{-(K+\alpha)t}\mathfrak{R}V(\phi(y,t))dt
\end{equation}
and
\begin{equation}
\label{defH}
H(V,y)=e^{-(K+\alpha)t^*(y)}\mathfrak{T}V(\phi(y,t^*(y))
\end{equation}
We define a time interval $\Delta t$ (in practice equal to one day) and for every $y=(\gamma,n,\sigma,\theta,p,r) \in \tilde{X}$, we note 
\begin{equation*}
n^*(y)=\left\lfloor \dfrac{t^*(y)}{\Delta t} \right\rfloor
\end{equation*}
For every $j \in \{0..n^*(y)-1\}$, we denote $\phi_j(y,t)=\phi(y,j\Delta t)$ and $\phi(y,t^*(y))=(\gamma,n,\sigma+t^*(y),\theta+t^*(y),p^*(y),r^*(y))$. The integral defined in equation \ref{defG} is computed by approximation using the classic trapezoidal rule using the $j\Delta t$ nodes. Thus, $G(V,y)$ can be approximated by a linear combination of $\{V(y_j)\}_{j \in \{0..n^*(y)-1\}}$, with $y_j$ depending on $\phi_j(y,t)$. Moreover, $H(V,y)$ is proportional to $V(\overline{y})$, with $\overline{y}$ depending on the boundary reached in $\phi(y,t^*(y))$. Finally, for every point $y \in \tilde{X}$, if we note $\overline{y}=y_{n^*(y)}$, $\mathfrak{B}V(y)$ can be computed as a linear combination of $\{V(y_j)\}_{j \in \{0..n^*(y)\}}$.

\section{Numerical aspects of the dynamic programing method}
\label{sec:3}
From theorem \ref{maintheorem} we know that we need to compute the sequence $\{W_q\}_{q\in \mathbb{N}}$ such that for $y\in \overline{\mathbf{X}}$, $W_0(y)=-K_{A}\mathbf{1}_{A_{\varepsilon_1}}(y)-(K_{A}+K_{B})\mathbf{1}_{A_{\varepsilon_1}^{c}}(y)$ and $W_{q+1}(y)=\mathfrak{B}W_q(y)$ for $q \in \mathbb{N}$. The sequence converges to a function $W$ defined on $\overline{\mathbf{X}}$ that allows the determination of the optimal cost and the optimal protocol of injections achieving that cost. However, computation has to be done on a grid $\Gamma$ of the state space: this grid must contain points of the form $(\gamma,n,\sigma,\theta,p,r)$. $\gamma$ and $n$ are discrete variables with $\gamma \in \{1...m_d+1\}$, $n \in \{1..n_{inj}\}$. $\sigma$ and $\theta$ are discretized with a time step of 1 day, with $\sigma \in \{0..\sigma_{\max}\}$ and $\theta \in \{0..T_h\}$. We then need to discretize $p$ and $r$: we choose a regular grid, with $p \in \{p_{\min}..p_{\max}\}$ with regular step $h_p$ and $r \in \{r_{\min}..r_{\max}\}$ with regular step $h_r$. We then obtain :
\begin{equation*}
n_p=\frac{p_{\max}-p_{\min}}{h_p}+1
\end{equation*}
\begin{equation*}
n_r=\frac{r_{\max}-r_{\min}}{h_r}+1
\end{equation*}
$h_r$ and $h_p$ are chosen such that both $n_p$, $n_r \in \mathbb{N}$ count the number of values of $p$ and $r$ on the grid, respectively. We arrange all points of the grid $\Gamma$ in a matrix $M$ of size $N_{sum} \times N_{pr}$, with $N_{sum}$ corresponding to the number of possible $(\gamma,n,\sigma,\theta)$ combinations and $N_{pr} = n_p n_r$ number of possible $(p,r)$ combinations. Each element $M(v,s) _{\begin{subarray}{l} v \in \{1..N_{sum}\}\\ s \in \{1..N_{pr}\} \end{subarray}}$ corresponds to a given combination $(\gamma,n,\sigma,\theta,p,r)$ of $\Gamma$, through the following bijection :
\begin{equation*}
    \begin{array}{ll}
         \chi \colon& \Gamma \to \{1..N_{sum}\} \times \{1..N_{pr}\}\\
        &(\gamma,n,\sigma,\theta,p,r) \mapsto (v,s)=\Big ( \chi_l(\gamma,n,\sigma,\theta), \chi_c(p,r)\Big )
    \end{array}
\end{equation*}
$\chi_l$ is defined in the following way : $v$ corresponds to a given value of $(\gamma,n,\sigma,\theta)$. Possible combinations of $(\sigma,\theta)$ depend on the value of $(\gamma,n)$: for example, during the first cycle, when $n=1$, $\sigma=0$ is associated with $\theta=1$, while when $n=2$, $\sigma=0$ is associated with $\theta=8$. We divide the lines of matrix $M$ by defining then $N_{\gamma n}= (m_d+1) n_{inj}$ blocks, corresponding to the possible combinations of $(\gamma,n)$. Each block is indexed by $i=f(\gamma,n)=\gamma + (m_d+1)(n-1) \in \{1..N_{\gamma n}\}$ and contains combinations of $(\sigma,\theta)$, indexed by $j=g_i(\sigma,\theta) \in \{1..N_{b_i}\}$ within the $i$-th block. Then $N_{sum} = \sum_{i=1}^{N_{\gamma n}} N_{b_i}$ and we can define a vector $l_{block}=(1,1+N_{b_1},..,1+\sum_{i=1}^{k} N_{b_i},.., \sum_{i=1}^{N_{\gamma n}-1} N_{b_i})$ of length $N_{\gamma n}$, that determines the index of the first line of each block. Finally : 
\begin{equation*}
v=\chi_l(\gamma,n,\sigma,\theta)=l_{block}(i)+j-1
\end{equation*}
with $i=f(\gamma,n)$ and $j=g_i(\sigma,\theta)$. $\chi_c$ is defined in the following way : 
\begin{equation*}
s=\chi_c(p,r)=\frac{p-p_{\min}}{h_p}+1+n_p \frac{r-r_{\min}}{h_r}
\end{equation*}
such that $s=1$ when $(p,r)=(p_{\min},r_{\min})$ and $s=n_p n_r$ when $(p,r)=(p_{\max},r_{\max})$.\\
Each iteration of the algorithm computes then a matrix $M_q$ such that \begin{equation*}
M_q(v,s)=W_q\Big(M(v,s)\Big )
\end{equation*}
For every $x=(\gamma,n,\sigma,\theta,p,r) \in \Gamma$, $W_q(x)$ is a linear combination of some $W_q(x_m)$, $m \in \{1..M_x\}$, as shown in equation \ref{defB} from appendix \ref{computeB}. Values of $W_q(x_m)$ are given by $M_q(\chi(x_m))$ and can then be combined and implemented in $M_{q+1}(\chi(x))$.

\section{Results}
\label{sec:4}
\begin{figure}
\centering
\includegraphics[width=1\textwidth]{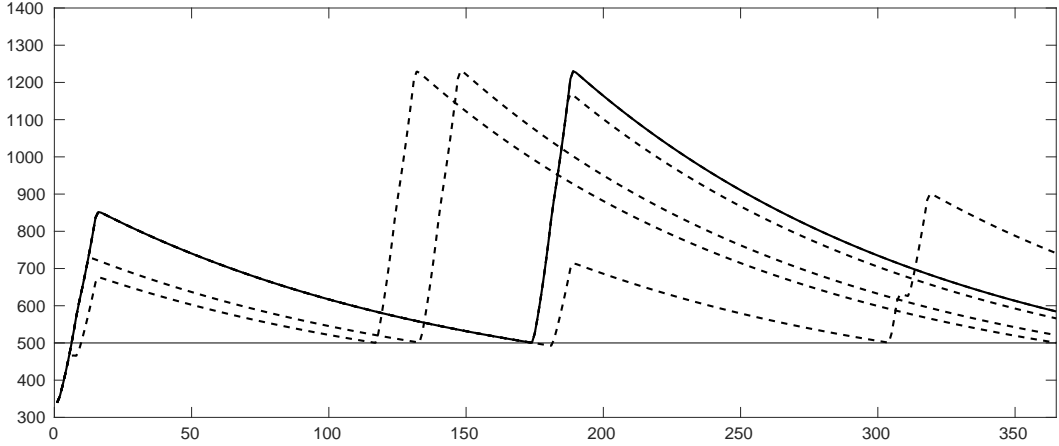}
\caption{Dynamics of CD4$^+$ T lymhocytes in patient A under a 2-injections cycles strategy (dose 20). Straight line correspond to the "best" outcome, i.e., when the effect of all injections lasts seven days. Dashed lines correspond to other possible trajectories, when this effect can last less than seven days}
\label{2injcycle}       
\end{figure}

\begin{figure}
\centering
\includegraphics[width=1\textwidth]{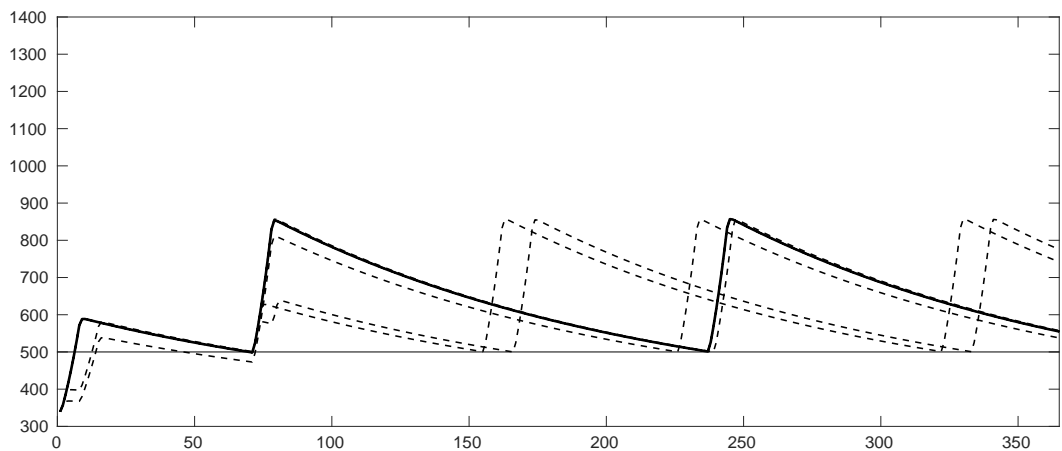}
\caption{Dynamics of CD4$^+$ T lymhocytes in patient A under the determined optimal strategy. Straight line correspond to the "best" outcome, i.e., when the effect of all injections lasts seven days. Dashed lines correspond to other possible trajectories, when this effect can last less than seven days}
\label{opt}       
\end{figure}

We applied the previously described method to a reduced model: as described in section \ref{sec:1.3}, we worked on a process with a choice of $m_d=2$ possibles doses: $d=[0,10,20]$ (unit = $\mu$g/kg), cycles of 2 injections: $n_{inj}=2$ and horizon $T_h=365$ days. We also assumed a minimum time of $\sigma_{\min}=70$ days between the beginning of two consecutive cycles. For a given patient with fixed biological parameters, we can approximate the function $W$ in a grid of the state space through convergence of the sequence $\{W_q\}$ : this determines the optimal cost over all strategies. Moreover, using equation \ref{action_eq} from theorem \ref{maintheorem}, we can simulate the strategy choosing the optimal action to realize when reaching the boundary of the state space and compute the cost of the obtained strategy. As some randomness is included in the model by the time of effect of an injection of parameter $\pi$, we assess the mean cost of a strategy by a Monte Carlo method using $N=10 000$ simulations of a protocol on a given patient. This allows the computation of the expectation of the cost of the strategy. From that, we check the performance of our method by first comparing the cost of the optimal strategy to the computation of the optimal cost from the value function $W$. Moreover, we wish to compare the optimal strategy to other "naive" protocols. For each protocol, including the optimal one, we compute by Monte Carlo the mean cost, the standard deviation and the minimum cost achievable. This is usually reached when the patient responds well to all injections, i.e. the effect of the injection on parameter $\pi$ lasts 7 days after every injection. In order to compare protocols based on clinical criteria, we also computed by Monte Carlo the mean number of CD4+ T cells count until horizon, the mean time spent under 500 cells/µL (in days) and the mean number of injections over all simulations. These comparisons were realized on two fictive patients, A and B with respectively average and low levels of CD4 at baseline, with values for parameters $\pi$, $\mu_R$, $\mu_P$ taken as mean value estimated in population in \cite{thiebaut2016repeated} and $\lambda$ and $\rho$ sampled from a normal law considering the variance of random effects previously estimated. Table \ref{patparam} sums up these two patients' characteristics. 
\begin{table}
\centering
\caption{Parameters for patients A and B}
\label{patparam}      
\begin{tabular}{lll}
\hline\noalign{\smallskip}
Parameter & Patient A  & Patient B \\
&(average baseline)&(low baseline)\\
\noalign{\smallskip}\hline\noalign{\smallskip}
$\lambda$ (cells/day) & 2.55 & 1.86 \\
$\rho$ (/day) & 2.06 & 1.10 \\
\hline
$\pi_0$ (/day) & 0.049 & 0.049 \\
$\mu_R$ (/day) & 0.054 & 0.054 \\
$\mu_P$ (/day) & 0.068 & 0.068 \\
\noalign{\smallskip}\hline
$\beta_{\pi_1}$ & 0.918 & 0.918 \\
$\beta_{\pi_2}$ & 0.721& 0.721 \\
\noalign{\smallskip}\hline
$R (t=0)$  &332 & 187 \\
$P(t=0)$ & 8  & 8 \\ 
\noalign{\smallskip}\hline
\end{tabular}
\end{table}
We first compare the value of the optimal function obtained from the numerical computation of $W$ with the cost of the optimal strategy. This result is displayed in table \ref{costcomp}. We note that for both patients the two cost values are very similar, meaning that we make a good approximation of the value function with our numerical method. 
\begin{table}
\centering
\caption{Comparison of cost values from value function and Monte Carlo simulation. Std = standard deviation}
\label{costcomp}      
\begin{tabular}{lll}
\hline\noalign{\smallskip}
& Patient A  & Patient B \\
\noalign{\smallskip}\hline\noalign{\smallskip}
 Optimal cost $W(x_0) $ & 3.22 & 9.65\\
\noalign{\smallskip}\hline
 Cost of optimal strategy : mean (std) &3.22 (0.56) & 9.66 (1.14)\\
 (obtained by Monte Carlo) & & \\
 \noalign{\smallskip}\hline
\end{tabular}
\end{table}
Then, we realized comparisons of several protocols. We simulated four "naive" protocols : one with cycles of two injections of dose 20, one with cycles of two injections of dose 10, one with cycles of one injection of dose 20 and finally one with a first cycle of two injections followed by cycles of one injection, all of dose 20. Assessing the cost of these protocols is interesting as they imply variable trajectories within the same patient as well as different values for clinical criteria. Moreover they would be clinically feasible and represent thus a good basis for comparison for our optimal strategy. Comparison results with cost function and clinical criteria are displayed respectively in table \ref{patmoyenAcf} and table \ref{patmoyenAcc} for patient A and in table \ref{patmoyenBcf} and \ref{patmoyenBcc} for patient B.
\begin{table}
\centering
\caption{Comparison of several protocols of injections for patient A using cost function. Min = minimum}
\label{patmoyenAcf}
\begin{tabular}{lccc}
\hline\noalign{\smallskip}\\
Protocol & Mean cost & Std & Min cost  \\
\hline\noalign{\smallskip}
2-injections cycles dose 20 & 4.15 &  0.70& 3.87\\
\hline\noalign{\smallskip}
2-injections cycles dose 10 &5.58  &0.50 & 5.42\\
\hline\noalign{\smallskip}
1-injection cycles dose 20 & 3.40 & 0.98 & 2.98 \\
\hline\noalign{\smallskip}
1st cycle : 2 injections & 3.87 & 0.45& 3.74  \\
then 1 injection, all doses 20&  &  &   \\
\hline\noalign{\smallskip}
\hline\noalign{\smallskip}
Optimal dynamical strategy : & 3.22 &  0.56& 2.98  \\
mostly 1-injection cycles dose 20. If response to 1st  & & & \\
injection is too low then 1st cycle of 2 injections&&&\\
\hline\noalign{\smallskip}
\end{tabular}
\end{table}
\begin{table}
\centering
\caption{Comparison of several protocols of injections for patient A using clinical criteria}
\label{patmoyenAcc}
\begin{tabular}{lccc}
\hline\noalign{\smallskip}
Protocol & CD4 mean & Days under 500 & Number of injections \\
\hline\noalign{\smallskip}
2-injections cycles dose 20 &  731 & 7.45&  4.29\\
\hline\noalign{\smallskip}
2-injections cycles dose 10 & 668& 10.9 &6.01 \\
\hline\noalign{\smallskip}
1-injection cycles dose 20 &  628 & 16.7& 3.18\\
\hline\noalign{\smallskip}
1st cycle : 2 injections & 646 &9.11  &4.02 \\
then 1 injection, all doses 20&  &  &   \\
\hline\noalign{\smallskip}
\hline\noalign{\smallskip}
Optimal dynamical strategy & 633 & 9.43 & 3.23 \\
\hline\noalign{\smallskip}
\end{tabular}
\end{table}
\begin{table}
\centering
\caption{Comparison of several protocols of injections for patient B using cost function}
\label{patmoyenBcf}
\begin{tabular}{lccc}
\hline\noalign{\smallskip}
Protocol & Mean cost & Std & Min cost \\
\hline\noalign{\smallskip}
2-injections cycles dose 20 & 10.48 & 1.18 & 9.64 \\
\hline\noalign{\smallskip}
2-injections cycles dose 10 & 13.69 & 1.04 & 13.01 \\
\hline\noalign{\smallskip}
1-injection cycles dose 20 & 13.78 & 0.63 &13.33 \\
\hline\noalign{\smallskip}
1st cycle : 2 injections & 13.06 & 1.19 & 12.18  \\
then 1 injection, all doses 20&  &  &  \\
\hline\noalign{\smallskip} 
\hline\noalign{\smallskip}
Optimal dynamical strategy :& 9.66 & 1.14 & 8.80  \\
2 cycles of 2 injections then 1-injection cycles.  &&& \\
If response too low then 3rd cycle of 2 injections &&&\\
\hline\noalign{\smallskip}
\end{tabular}
\end{table}
\begin{table}
\centering
\caption{Comparison of several protocols of injections for patient B using clinical criteria}
\label{patmoyenBcc}
\begin{tabular}{lccc}
\hline\noalign{\smallskip}
Protocol  & CD4 mean & Days under 500 & Number of injections \\
\hline\noalign{\smallskip}
2-injections cycles dose 20 & 620 & 92.5 & 8.54 \\
\hline\noalign{\smallskip}
2-injections cycles dose 10 & 522 & 162 & 10.1\\
\hline\noalign{\smallskip}
1-injection cycles dose 20  & 398 & 306 & 6 \\
\hline\noalign{\smallskip}
1st cycle : 2 injections & 462 & 242 & 7\\
then 1 injection, all doses 20&  &  &  \\
\hline\noalign{\smallskip} 
\hline\noalign{\smallskip}
Optimal dynamical strategy & 561 & 92.9 & 7.57 \\
\hline\noalign{\smallskip}
\end{tabular}
\end{table}
A first result that should be underlined is that for both patients, mean cost of the optimal strategy is lower than all other simulated strategies. We also observe that the cost is higher in all strategies for patient B than for patient A, which is consistent with the fact that patient B has baseline CD4 levels lower than patient A and thus will spend more time under 500 cells/ $\mu$L and need more injections to raise above it. For patient A, the optimal strategy is achieved by one-injection cycles, except when patient's response to the first injection is too low. In that case, the first cycle contains two injections and the next cycles includes only one injection. This explains why the minimum cost is the same for the optimal strategy and the one-injection cycles at dose 20 strategy (that is equal to 2.98): when the patient has a good response to all injections, these are the same strategies. Of note, for patient A the optimal strategy achieves a good balance between all clinical criteria, even if this strategy does not seem optimal when the criteria are taken separately. This is a good point, as the cost function was chosen in order to minimize both criteria. In the case of patient B, the optimal strategy contains two first cycles of two injections followed by cycles of one injection, and if the number of CD4 are still too low after those two cycles, the third one has two injections as well (all doses = 20  $\mu$g/kg). This strategy is very intuitive : the first complete cycles are needed to raise the number of CD4 over 500 cells/ $\mu$L, as the baseline levels are very low (less than 200 cells/ $\mu$L); then, one-injection cycles allow to sustain the levels over 500 cells/ $\mu$L. We can note that, as in the case of patient A, the optimal strategy allows a good balance between the number of injections realized and the time spent under 500 cells/ $\mu$L. For example, the optimal strategy allows to spend as much time over 500 cells/ $\mu$L as the 2-injections cycles (dose 20  $\mu$g/kg) strategy by using one less injection. All together, these results suggest that our numerical method allows to determine an optimal strategy of injections, and the clinical interpretation of the results are consistent with the mathematical method. In term of trajectories of the process, figures \ref{2injcycle} and \ref{opt} show some trajectories obtained with respectively the 2-injection cycles strategy and the optimal one for patient A. We can note that even if CD4$^+$ levels are globally lower in the optimal strategy compared to the two-injections cycles at dose 20  $\mu$g/kg, it still allows a maintenance over the threshold of 500 cells/$\mu$L by using less injections: indeed, in the best case scenario, the 2-injections cycle strategy implies 2 cycles of 2 injections which is a total of 4 injections, while the optimal strategy induces 3 cycles of 1 injection, which is a total of 3 injections. This supports the interest of determining the optimal strategy based on a criterion combining both the number of injections and the time spent under 500.

\section{Discussion}
\label{sec:5}
In this work, we have developed a numerical tool allowing to solve an impulse control problem for a PDMP. The specificity of our work is in the application of the dynamic programming approach to a specific biological framework. The objective is to determine the optimal strategy of IL-7 injections for a given HIV-infected patient, in order to maintain CD4$^+$ T lymphocytes levels over the threshold of 500 cells/$\mu$L. We first modeled the dynamics of CD4$^+$ T lymphocytes during repeated cycles of IL-7 injections by a PDMP. Then, we solved the impulse control problem by iterating a sequence defined by an integro-differential operator. This sequence converges to the value function, which allows to determine the optimal action that should be realized at every point of the boundary. Thus, we propose a numerical tool approximating the sequence and the value function on a grid of the state space. We applied this tool to a simplified clinical question. Although the horizon of study is only one year and the number of injections limited to two by cycle, we obtained interesting results, also consistent with a clinical interpretation. The optimal strategy determined for two different patients is indeed intuitive: for both of them, the first cycles aim at increasing the CD4$^+$ T lymphocytes levels and should contain as many injections as possible until the levels are acceptable. Then, the following cycles sustain the CD4 levels over the threshold, and punctual injections are sufficient to reach this objective. The optimal strategy, determined with our method, has a lower cost than other possible clinical strategies. Actually, the obtained optimal strategy depends on the cost previously defined, and thus we could explore other optimal strategies depending on other cost functions. For example, it could be interesting to use different weights on the time spent under 500  cells /$\mu$L and the number of injections (depending on the clinician priorities), or to account for the possible negative side effects due to higher doses (this would need additional data on the question). Finally, the model could be extended, accounting for three possible injections by cycle (and not only two) and studying the patient until a longer horizon (up to two years). This rises the issue of the increase of computational time (by increasing the size of the grid of the state space) and constitutes a new challenge in itself. In the end, we hope to use this tool in future possible clinical trial investigating the effect of IL-7 injections with patients-specific schedules of injections, personalized and optimized using this method.

\section*{Acknowledgments}

We would like to thank the main investigators and supervisors of INSPIRE 2 and 3 studies: Jean-Pierre Routy, Irini Sereti,
Margaret Fischl, Prudence Ive, Roberto F. Speck, Gianpiero D'Ozi, Salvatore Casari, Sharne Foulkes, Ven Natarajan, Guiseppe Tambussi, Michael M. Lederman, Therese Croughs and Jean-Fran\c{c}ois Delfraissy. This work was supported by the Investissements d'Avenir program managed by the ANR under reference ANR-10-LABX-77.

\bibliographystyle{spbasic}      
\bibliography{bibcontrole}   

%
%

\appendix
\section{Optimal control : application}
\label{computeB}
We defined the process describing the effect of IL-7 on CD4$^+$ T lymphocytes dynamics by its characteristics $\phi$, $\eta$ and $Q$, boundaries and possible actions in section \ref{sec:1.3}. We also defined both gradual cost on the trajectory and impulse cost in that section. As we aim at applying the results from theorem \ref{maintheorem} to determine the optimal cost and an optimal strategy by dynamic programming, we need to determine how to compute numerically the function $\mathfrak{B}$ to iterate the sequence $\{W_q\}_{q \in \mathbb{N}}$ defined in equation \ref{seqW}. As a reminder, $\mathfrak{B}$ is defined in \cite{costa2016constrained} by :
\begin{equation*}
\mathfrak{B}V(y)=
\int_{[0,t^{*}(y)[} e^{-(K+\alpha)t} \mathfrak{R}V(\phi(y,t)) dt + e^{-(K+\alpha)t^{*}(y)} \mathfrak{T}V(\phi(y,t^{*}(y)))
\end{equation*}
We will first detail the computation of $\mathfrak{R}$ then $\mathfrak{T}$, and we will finally show how to compute $\mathfrak{B}$. 
\subsection*{Computation of $\mathfrak{R}$}
For $x=(\gamma,n,\sigma,\theta,p,r) \in X $, and function $V : \overline{X} \rightarrow \mathbb{R}$, $\mathfrak{R}$ is defined as :
\begin{equation*}
\mathfrak{R}V(x) = C^g(x) + qV(x) + \eta V(x)   
\end{equation*}
with q computing the difference between the states before and after the spontaneous jump. As Q depends on the action only when the process hits the active boundary, 
\begin{equation*}
\begin{array}{rll}
        q(dy|x,d)& =&\eta(x)[Q(dy|x)-\delta_x(dy)]\\
       &=&\mathbf{1}_{\{\gamma>1\}}\eta[\delta_{(1,n,\sigma,\theta,p,r)}(dy)-\delta_{(\gamma,n,\sigma,\theta,p,r)}(dy)]
\end{array}
\end{equation*}
then for every function $V $, and as $K=\eta$ : 
\begin{equation*}
\begin{array}{rll}
qV(x) &=& \displaystyle\int V(y)q(dy|x) \nonumber\\
&=& \mathbf{1}_{\{\gamma>1\}}K[V(1,n,\sigma,\theta,p,r)-V(\gamma,n,\sigma,\theta,p,r)]
\end{array}
\end{equation*}
Then 
\begin{equation*}
\begin{array}{rll}
\mathfrak{R}V(x)&=& \displaystyle\frac{1}{30}\mathbf{1}_{\{p+r\leq 500\}} + qV(x) + K V(x)  \\
&= & \displaystyle\frac{1}{30}\mathbf{1}_{\{p+r\leq 500\}} +K V(1,n,\sigma,\theta,p,r)  \mathbf{1}_{\{\gamma>1\}} + K V(x)  \mathbf{1}_{\{\gamma=1\}}  \\
\end{array}
\end{equation*}
Finally,
\begin{equation}
\label{Rdef}
\begin{array}{rll}
\mathfrak{R}V(x) &=& \displaystyle\frac{1}{30} \mathbf{1}_{\{p+r\leq 500\}} + K V(1,n,\sigma,\theta,p,r) \\
\mathfrak{R}V(\Delta) &=& KV(\Delta)
\end{array}
\end{equation}
\subsection*{Computation of $\mathfrak{T}$}
For $x \in \Xi $, and function $V : \overline{X} \rightarrow \mathbb{R}$, $\mathfrak{T}$ is defined as : 
\begin{equation*}
\begin{array}{rll}
\mathfrak{T}V(x) &=& \displaystyle\inf_{d \in A(x)} \Big\{C^i(x,d) + QV(x,d) \Big\} \\
&=&\displaystyle\inf_{d \in A(x)}  \Big\{\mathbf{1}_{x \in \Xi_1 \cup \Xi_4} + \mathbf{1}_{d \neq 0}\mathbf{}{1}_{x \in \Xi_3} + \int V(y)\Big [\delta_{(\gamma(d),1,0,1,P_c,R_c)}(dy)\mathbf{1}_{\{x \in \Xi_1\}} \\
&&+ \delta_{\Delta}(dy)\mathbf{1}_{\{x \in \Xi_2\}} 
 + \delta_{(\gamma(d),n+1,0,\theta,p,r)}(dy)\mathbf{1}_{\{x \in \Xi_3\}} + \delta_{(\gamma(d),1,0,\theta,p,r)}(dy)\mathbf{1}_{\{x \in \Xi_4\}} \\
 &&+\delta_{(1,n,\sigma,\theta,p,r)}(dy)\mathbf{1}_{\{x \in \Xi_5\}}\Big ] \Big\} 
\end{array}
\end{equation*}
Finally, 
\begin{equation}
\label{defT}
\begin{array}{rll}
 \mathfrak{T}V(x) &=& \displaystyle\inf_{d \in A(x)}  \Big\{ [1+V(\gamma(d),1,0,1,P_c,R_c)] \mathbf{1}_{x \in \Xi_1} \\
 &&+ [\mathbf{1}_{d \neq 0}+V(\gamma(d),n+1,0,\theta,p,r)]  \mathbf{1}_{x \in \Xi_3}  \\
 & &+ [1+V(\gamma(d),1,0,\theta,p,r)]  \mathbf{1}_{x \in \Xi_4} \Big\} + V(\Delta) \mathbf{1}_{x \in \Xi_2} + V(1,n,\sigma,\theta,p,r) \mathbf{1}_{x \in \Xi_5} \\
 \mathfrak{T}V(\Delta)&=&V(\Delta)
\end{array}
\end{equation}
\subsection*{Computation of $\mathfrak{B}$}
Now, for $Y \in \overline{X} $, and function $V : \overline{X} \rightarrow \mathbb{R}$, we need to compute :
\begin{equation*}
\mathfrak{B}V(y)=
\int_{[0,t^{*}(y)[} e^{-(K+\alpha)t} \mathfrak{R}V(\phi(y,t)) dt + e^{-(K+\alpha)t^{*}(y)} \mathfrak{T}V(\phi(y,t^{*}(y)))
\end{equation*}
As we cannot make an exact computation of $\mathfrak{B}V$ on $\overline{X}$, we need to approximate this computation on a grid of the state space. In order to detail the approximation of the computation, we define 
\begin{equation*}
G(V,y)=\int_{[0,t^*(y)[} e^{-(K+\alpha)t}\mathfrak{R}V(\phi(y,t))dt
\end{equation*}
and
\begin{equation*}
H(V,y)=e^{-(K+\alpha)t^*(y)}\mathfrak{T}V(\phi(y,t^*(y))
\end{equation*}
as in \ref{defG} and \ref{defH}. We define a time interval $\Delta t$ (in practice equal to one day) and for every $y=(\gamma,n,\sigma,\theta,p,r) \in \tilde{X}$, we note 
\begin{equation*}
n^*(y)=\left\lfloor \dfrac{t^*(y)}{\Delta t} \right\rfloor
\end{equation*}
For every $j \in \{0..n^*(y)-1\}$, we note $\phi_j(y,t)=\phi(y,j\Delta t)$ and $\phi(y,t^*(y))=(\gamma,n,\sigma+t^*(y),\theta+t^*(y),p^*(y),r^*(y))$. The integral defined in equation \ref{defG} is computed by approximation using the classic trapezoidal rule using the $j\Delta t$ nodes : 
\begin{equation*}
\begin{array}{rll}
G(V,y) 
&\simeq &\displaystyle\frac{\Delta t}{2}\mathfrak{R}V(y)+\frac{\Delta t}{2}e^{-(K+\alpha)t^*(y)}\mathfrak{R}V(\phi(y,t^*(y))) \\
&&+ \displaystyle\sum_{j=1}^{n^*(y)-2} \Delta t e^{-(K+\alpha)j \Delta t} \mathfrak{R}V(\phi_j(y,t))
\end{array}
\end{equation*}
with $\mathfrak{R}V(x) = \displaystyle\frac{1}{30} \mathbf{1}_{\{p+r\leq 500\}} + K V(1,n,\sigma,\theta,p,r)$, as computed in \ref{Rdef}. Then we obtain the following for every $y=(\gamma,n,\sigma,\theta,p,r) \in \tilde{X}$ : 
\begin{equation*}
\begin{array}{rll}
G(V,y)&=&\displaystyle\frac{\Delta t}{2}(\displaystyle\frac{1}{30} \mathbf{1}_{\{p+r<500\}} + K V(1,n,\sigma,\theta,p,r)) \\
&+& \displaystyle\frac{\Delta t}{2}e^{-(K+\alpha)t^*(y)}(\displaystyle\frac{1}{30}\mathbf{1}_{\{p^*+r^*<500\}}+ K V(1,n,\sigma +t^*,\theta +t^*,p^*(y),r^*(y)))  \\
&+& \Delta t \displaystyle\sum_{j=1}^{n^*(y)-2}  e^{-(K+\alpha)j \Delta t} (\displaystyle\frac{1}{30}\mathbf{1}_{\{p_j+r_j<500\}}+K V(1,n,\sigma + j\Delta t,\theta +j\Delta t,p_j,r_j)) 
\end{array}
\end{equation*}
Now, we need to compute H as defined in \ref{defH} : it depends on $\mathfrak{T}V(\phi(y,t^*(y)))$, which takes different values according to the boundary reached in that point, as written in equation \ref{defT}. Moreover, as we know the flow, we can give conditions on $y=(\gamma,n,\sigma,\theta,p,r)$ to reach a given boundary in $\phi(y,t^*(y))$. Then : 
\begin{itemize}
\item if $\phi(y,t^*(y)) \in \Xi_1$ ($\theta \leq 1$) then 
\begin{equation}
H(V,y)=  \inf_{d \in [d_1,..d_{m_d}]}  \Big\{ e^{-(K+\alpha)t^*(y)}\Big[1+V(\gamma(d),1,0,1,P_{c},R_{c}) \Big]\Big\}
\end{equation}
\item if $\phi(y,t^*(y)) \in \Xi_2$ ($\theta + t^*(y) \geq T_h$) then 
\begin{equation*}
H(V,y)= e^{-(K+\alpha)t^*(y)} V(\Delta)
\end{equation*}
\item if $\phi(y,t^*(y)) \in \Xi_3$ ($n < n_{inj}, \theta + t^*(y) < T_h$) then 
\begin{equation*}
\begin{array}{rll}
H(V,y)&=&  \inf_{d \in [0,d_1,..d_{m_d}]}  \Big\{ e^{-(K+\alpha)t^*(y)}\Big[\mathbf{1}_{\{d \neq 0\}}\\ 
&&+V(\gamma(d),n+1,0,\theta+t^*(y),p^*(y),r^*(y)) \Big]\Big\}
\end{array}
\end{equation*}
\item if $\phi(y,t^*(y)) \in \Xi_4$ ($n = n_{inj}, \gamma=1, \theta + t^*(y) < T_h$) then 
\begin{equation*}
H(V,y)=  \inf_{d \in [d_1,..d_{m_d}]}  \Big\{ e^{-(K+\alpha)t^*(y)}\Big[1+V(\gamma(d),1,0,\theta+t^*(y),p^*(y),r^*(y)) \Big]\Big\}
\end{equation*}
\item if $\phi(y,t^*(y)) \in \Xi_5$ ($n = n_{inj}, \gamma>1, \theta + t^*(y) < T_h$) then 
\begin{equation*}
H(V,y)= e^{-(K+\alpha)t^*(y)} V(1,n,\sigma+t^*(y),\theta+t^*(y),p^*(y),r^*(y))
\end{equation*}

\end{itemize}
Finally, for every $y=(\gamma,n,\sigma,\theta,p,r)\in \tilde{X}$ : 
\begin{equation}
\begin{array}{rl}
\label{defB}
\mathfrak{B}V&(y)=\displaystyle\frac{\Delta t}{2}(\displaystyle\frac{1}{30}\mathbf{1}_{\{p+r<500\}} + K V(1,n,\sigma,\theta,p,r))\\
& +\frac{\Delta t}{2}e^{-(K+\alpha)t^*(y)}\Big[\displaystyle\frac{1}{30}\mathbf{1}_{\{p^*+r^*<500\}}+ K V(1,n,\sigma +t^*,\theta +t^*,p^*(y),r^*(y))\Big]  \\
&+\Delta t \displaystyle\sum_{j=1}^{n^*(y)-2}  e^{-(K+\alpha)j \Delta t} \Big[\displaystyle\frac{1}{30}\mathbf{1}_{\{p_j+r_j<500\}}+K V(1,n,\sigma + j\Delta t,\theta +j\Delta t,p_j,r_j)\Big] \\ 
&+\displaystyle\inf_{d \in [d_1,..d_{m_d}]}  \Big\{ e^{-(K+\alpha)t^*(y)}\Big[1+V(\gamma(d),1,0,1,P_{c},R_{c}) \Big]\Big\} \mathbf{1}_{\{ \theta \leq 1 \}} \\
&+e^{-(K+\alpha)t^*(y)} V(\Delta) \mathbf{1}_{\{ \theta + t^*(y) \geq T_h\} }\\ 
&+ \displaystyle\inf_{d \in [0,d_1,..d_{m_d}]}  \Big\{ e^{-(K+\alpha)t^*(y)}\Big[\mathbf{1}_{\{d \neq 0\}}\\
&+V(\gamma(d),n+1,0,\theta+t^*(y),p^*(y),r^*(y)) \Big]\Big\} \mathbf{1}_{\{ n < n_{inj}, \theta + t^*(y) < T_h\}} \\ 
&+\displaystyle\inf_{d \in [d_1,..d_{m_d}]}  \Big\{ e^{-(K+\alpha)t^*(y)}\Big[1\\
&+V(\gamma(d),1,0,\theta+t^*(y),p^*(y),r^*(y)) \Big]\Big\} \mathbf{1}_{\{ n = n_{inj}, \gamma=1, \theta + t^*(y) < T_h\}} \\ 
&+ e^{-(K+\alpha)t^*(y)} V(1,n,\sigma+t^*(y),\theta+t^*(y),p^*(y),r^*(y)) \mathbf{1}_{\{n = n_{inj}, \gamma>1, \theta + t^*(y) < T_h\}}
\end{array}
\end{equation}
and 
\begin{equation*}
\mathfrak{B}V(\Delta)=\int_{[0,+\infty)} e^{-(K+\alpha)t}KV(\Delta)dt = \frac{K}{K+\alpha} V(\Delta)
\end{equation*}
\end{document}